\documentclass{tac}
\usepackage[square]{natbib}
\usepackage{amssymb}
\usepackage{graphicx}
\usepackage[all,2cell]{xy}
\UseAllTwocells
\xyoption{v2}
\NoResizing
\UseResizing
\usepackage[colorlinks=true]{hyperref}
\hypersetup{allcolors=[rgb]{0.1,0.1,0.4}}

\title{The Existential Completion}
\author{Davide Trotta}
\copyrightyear{2020}
\address{Department of Mathematics, University of Trento\\
Via Sommarive, 14, Trento, Italy}
\eaddress{trottadavide92@gmail.com}
\keywords{existential completion, doctrines, property-like monad, tripos}
\amsclass{18C10, 18C15, 18D05}

\thanks{My acknowledgements go first to my supervisor Pino Rosolini, for his indispensable comments and suggestions. Then, I wish also to thank Milly Maietti for discussions about examples and applications of the existential completion, and the anonymous referee for helpful feedback.}

\newcommand{\operatorname}[1]{\mathop{#1}\nolimits}
\newcommand{\bemph}[1]{\textbf{\emph{#1}}}
\newcommand{\alg}[1]{\operatorname{#1\mbox{-}\mathbf{Alg}}}
\newcommand{\alglax}[1]{\operatorname{#1\mbox{-}\textbf{Alg}_l}}
\newcommand{\compex}[1]{{#1}^e}
\newcommand{\ovln}[1]{\overline{#1}}
\newcommand{\angbr}[2]{\langle #1,#2 \rangle} 
\newcommand{\freccia}[3]{\ensuremath{\xymatrix{#2 \colon #1  \ar[r] &  #3}}}
\newcommand{\frecciasopra}[3]{\xymatrix{ #1  \ar[r]^{#2} &  #3}}
\newcommand{\frecciasopralunga}[3]{\xymatrix{ #1  \ar[rr]^{#2} &&  #3}}
\newcommand{\pbmorph}[2]{#2^{\ast}#1} 
\newcommand{\duefreccia}[3]{\xymatrix@C=0.5cm{#2 \colon #1  \ar@{=>}[r] &  #3}}
\newcommand{\doctrine}[2]{\xymatrix{#2 \colon #1^{op}  \ar[r] & \infsl }}
\newcommand{\duemorfismo}[6]{\xymatrix{
#1^{op} \ar[rrd]^#2_{}="a" \ar[dd]_{#3^{op}}\\
&& \infsl\\
#5^{op}  \ar[rru]_#6^{}="b"
\ar_#4  "a";"b"}}
\newcommand{\comsquare}[8]{ \xymatrix@+1pc{ 
#1 \ar[r]^{#5} \ar[d]_{#6} & #2 \ar[d]^{#7} \\
#3 \ar[r]_{#8} & #4 
}}
\newcommand{\quadratocomm}[8]{ \xymatrix{ 
#1 \ar[r]^{#5} \ar[d]_{#6} & #2 \ar[d]^{#7} \\
#3 \ar[r]_{#8} & #4 
}}
\newcommand{\pullbackcorner}[1][ul]{\save*!/#1+1.2pc/#1:(1,-1)@^{|-}\restore}
\newcommand{\pbp}[3]{#1\times_{#3} #2}
\newcommand{\TC}[1]{\ensuremath{\mathbf{#1}}\xspace}
\newcommand{\ED}{\TC{ED}}
\newcommand{\PD}{\TC{PD}}
\newcommand{\CEED}{\TC{CEED}}
\newcommand{\ElD}{\TC{ElD}}
\newcommand{\PED}{\TC{ED}}
\newcommand{\Cat}{\TC{Cat}}
\newcommand{\EED}{\TC{EED}}
\newcommand{\LFS}{\TC{LFS}}
\newcommand{\excat}{\TC{Xct}}
\newcommand{\SL}[1]{\ensuremath{\mathcal{#1}}\xspace}
\newcommand{\mA}{\SL{A}}
\newcommand{\mC}{\SL{C}}
\newcommand{\mD}{\SL{D}}
\newcommand{\mS}{\SL{H}}
\newcommand{\mT}{\mathrm{T}}
\newcommand{\forg}{\mathrm{U}}
\newcommand{\funE}{\mathrm{E}}
\newcommand{\funU}{\mathrm{U}}
\newbox\erove \setbox\erove=\hbox{\reflectbox{E}}
\newcommand{\Einv}{\usebox\erove}
\newcommand{\pr}{\operatorname{\hspace{0.01cm} pr}}
\newcommand{\id}{\operatorname{\hspace{0.01cm} id}}
\newcommand{\op}{\operatorname{\hspace{0.01cm} op}}
\newcommand{\Sub}{\operatorname{\mathsf{Sub}}}
\newcommand{\Reg}{\operatorname{\hspace{0.01cm} \textbf{Reg}}}
\newcommand{\theory}{\operatorname{\hspace{0.01cm} \mathcal{T}}}
\newcommand{\infsl}{\operatorname{\hspace{0.01cm} \mathbf{InfSL}}}
\newcommand{\pos}{\operatorname{\hspace{0.01cm} \mathbf{Pos}_{\top}}}
\newcommand{\lang}{\operatorname{\mathcal{L}}}
\begin{document}
\maketitle
\begin{abstract}
We determine the existential completion of a primary doctrine, and we
prove that the 2-monad obtained from it is lax-idempotent,
and that the 2-category of existential doctrines is isomorphic to the
2-category of algebras for this 2-monad. We also show that the
existential completion of an elementary doctrine is again elementary.
Finally we extend the notion of exact completion of an elementary
existential doctrine to an arbitrary elementary doctrine.

\end{abstract}
\section{Introduction}
In recent years, many relevant logical completions have been
extensively studied in category theory. The main instance is the
exact completion, see \citep{SFEC,FECLEO,REC}, which is the universal
extension of a category with finite limits to an exact category.
In \citep{EQC,QCFF,UEC}, Maietti and Rosolini introduce a categorical
version of quotient for an equivalence relation, and they study that
in a doctrine equipped with a sufficient logical structure to describe
the notion of an equivalence relation.
In \citep{UEC} they show that both the exact completion of a regular
category and the exact completion of a category with binary products,
a weak terminal object and weak pullbacks can be seen as instances of
a more general completion with respect to an elementary existential
doctrine.

In this paper we present the existential completion of a primary
doctrine, and we give an explicit description of the 2-monad
$\freccia{\PD}{\mT_e}{\PD}$ constructed from the 2-adjunction, where
$\PD$ is the 2-category of primary doctrines.

It is well known that pseudo-monads can express uniformly and
elegantly many algebraic structures; we refer the reader to
\citep{UCTSBSSL,PDLAVB,OPLS} for a detailed description of these
topics. 

Recall that an action of a 2-monad on a given object encodes a \emph{structure} on that object. When the structure is uniquely determined to within unique isomorphism, to give an object with such a structure is just to give an object with a certain \emph{property}.
Those 2-monads for which the algebra structure is essentially unique, if it exists, are called \emph{property-like}.

In this paper we show that every existential doctrine $\doctrine{\mC}{P}$
admits an action $\freccia{\mT_e(P)}{a}{P}$ such that $(P,a)$ is a
$\mT_e$-algebra, and that if $(R,b)$ is $\mT_e$-algebra then the
doctrine is existential, and this gives an equivalence between the
2-category $\alg{\mT_e}$ and the 2-category $\ED$ whose objects are
existential doctrines.

Here the action encodes the existential structure for a doctrine,
and we prove that this structure is uniquely determined to within
appropriate isomorphism, i.e. that the 2-monad $\mT_e$ is lax-idempotent and hence property-like in the sense of \citep{OPLS}.

We also prove that the existential completion
preserves the elementary structure of a doctrine, and then we generalize the
bi-adjunction $\EED\rightarrow \excat$ presented in \citep{UEC,TECH} to
a bi-adjunction from the 2-category $\ElD$ of 
elementary doctrines to the 2-category of exact categories $\excat$.

In the sections 2 and 3 we recall definitions and results
on 2-monads, and on primary and existential
doctrines as needed for the rest of the paper.

In section 4 we describe the existential completion.
We introduce a 2-functor from the
2-category of primary doctrines to the 2-category of existential
doctrines $\freccia{\PD}{\funE}{\PED,}$ and we prove that it is a left 2-adjoint 
to the forgetful functor $\freccia{\PED}{\funU}{\PD.}$

In section 5 we prove that the 2-monad $\mT_e$ constructed
from the 2-adjunction is lax-idempotent and that the
2-category $\alg{\mT_e}$ is 2-equivalent to the 2-category $\PED$ of
existential doctrines.

In section 6 we show that the existential completion preserves the elementary structure, and we use this result to extend the notion of
exact completion to elementary doctrines.

\section{A brief recap of two-dimensional monad theory}
This section is devoted to the formal definition of 2-monad on a
2-category and a characterization of the definitions. We use
2-categorical pasting notation freely, following the usual convention
of the topic as used extensively in \citep{TDMT}, \citep{PDLAVB} and
\citep{UCTSBSSL}.

You can find all the details of the main results of this section in
the works of Kelly and Lack \citep{OPLS}. For a more general and
complete description of these topics, and a generalization for the
case of pseudo-monad, you can see the Ph.D thesis of Tanaka
\citep{PHDTT}, the articles of Marmolejo \citep{CPLR}, \citep{DLP} and
the work of Kelly \citep{RE2C}. Moreover we refer to \citep{HCA1} and
\citep{HOHC} for all the standard results and notions about 2-category
theory.

A \bemph{2-monad} $(\mT,\mu, \eta)$ on a 2-category $\mA$ is a 2-functor $\freccia{\mA}{\mT}{\mA}$ together 2-natural transformations $\freccia{\mT^2}{\mu}{\mT}$ and $\freccia{1_{\mA}}{\eta}{\mT}$ such that the following diagrams 
\[\xymatrix@+1pc{
\mT^3\ar[d]_{\mu\mT} \ar[r]^{\mT\mu} &\mT^2\ar[d]^{\mu}\\
\mT^2 \ar[r]_{\mu} & \mT
}\]
\[\xymatrix@+1pc{
\mT \ar[dr]_{\id}\ar[r]^{\eta \mT}& \mT^{2} \ar[d]_{\mu} & \mT \ar[l]_{\mT \eta}\ar[dl]^{\id}\\
& \mT
}\]
commute.

Let $(\mT,\mu, \eta)$ be a 2-monad on a 2-category $\mA$. A $\mT$-\bemph{algebra} is a pair $(A,a)$ where, $A$ is an object of $\mA$ and $\freccia{\mT A}{a}{A}$ is a 1-cell such that the following diagrams commute
\[\xymatrix@+1pc{
\mT^2 A  \ar[r]^{\mT a} \ar[d]_{\mu_A} & \mT A \ar[d]^{a}\\
\mT A \ar[r]_a & A
}\]
\[\xymatrix@+1pc{
A \ar[rd]_{1_A} \ar[r]^{\eta_A}& \mT A \ar[d]^{a}\\
& A.
}\]

A \bemph{lax} $\mT$-\bemph{morphism} from a $\mT$-algebra $(A,a)$ to a $\mT$-algebra $(B,b)$ is a pair $(f,\ovln{f})$ where $f$ is a 1-cell $\freccia{A}{f}{B}$ and $\ovln{f}$ is a 2-cell
\[\xymatrix@+1pc{
\mT A \ar[d]_{a} \ar[r]^{\mT f} \xtwocell[r]{}<>{_<5>\ovln{f}}& \mT B \ar[d]^b\\
A \ar[r]_f &B
}\] 
which satisfies the following \bemph{coherence} conditions
\[\xymatrix@+1pc{
\mT^2A \ar[d]_{\mu_A}\ar[r]^{\mT^2f} &TB \ar[d]^{\mu_B}&&\mT^2A \ar[d]_{\mT a}\xtwocell[r]{}<>{_<5>\;\;\;\mT\ovln{f}}\ar[r]^{\mT^2f} &\mT B \ar[d]^{\mT b}\\ 
\mT A \ar[d]_{a} \ar[r]^{\mT f} \xtwocell[r]{}<>{_<5>\ovln{f}}&\mT B \ar[d]^b\ar@{}[rr]|{=}&& \mT A \ar[d]_{a} \ar[r]^{\mT f} \xtwocell[r]{}<>{_<5>\ovln{f}}&\mT B \ar[d]^b \\
A \ar[r]_f &B&& A \ar[r]_f &B
}\]
and
\[\xymatrix@+1pc{
A \ar[d]_{\eta_A}\ar[r]^f &B\ar[d]^{\eta_B}&& A\ar[dd]_{1_A}\ar[r]^f& B\ar[dd]^{1_B}\\
\mT A \ar[d]_{a} \ar[r]^{\mT f} \xtwocell[r]{}<>{_<5>\ovln{f}}&\mT B \ar[d]^b\ar@{}[rr]|{=}&&\\
A \ar[r]_f & B&& A\ar[r]_f &B.
}\]
The regions in which no 2-cell is written always commute by the
naturality of $\eta$ and $\mu$, and are deemed to contain the identity
2-cell.

A lax morphism $(f,\ovln{f})$ in which $\ovln{f}$ is invertible
is said $\mT$-\bemph{morphism}. And it is \bemph{strict} when
$\ovln{f}$ is the identity.

The category of $\mT$-algebras and lax $\mT$-morphisms becomes a
2-category $\alglax{\mT}$, when provided with 2-cells the
$\mT$-\bemph{transformations}. Recall from \citep{OPLS} that
a $\mT$-\bemph{transformation} 
from $\freccia{(A,a)}{(f,\ovln{f})}{(B,b)}$ to
$\freccia{(A,a)}{(g,\ovln{g})}{(B,b)}$ is a 2-cell 
$\duefreccia{f}{\alpha}{g}$ in $\mA$ which satisfies the following
coherence condition 
\[\xymatrix{
\mT A \ar[dd]_{a}\rrtwocell<4>^{\mT f}_{\mT g}{\;\;\;\mT\alpha} \xtwocell[rr]{}<>{_<7>\ovln{g}}& &\mT B\ar[dd]^b&&\mT A \ar[dd]_{a}\ar@/^/[rr]^{\mT f} \xtwocell[rr]{}<>{_<7>\ovln{f}}& &\mT B\ar[dd]^b\\
&&\ar@{}[rr]|{=}&&\\
A \ar@/_/[rr]_g&& B&& A\rrtwocell<4>^f_g{\alpha} && B
}\]
expressing compatibility of $\alpha$ with $\ovln{f}$ and $\ovln{g}$.

It is observed in \citep{OPLS} that using this notion of
$\mT$-morphism, one can express more precisely what it may mean that
an action of a monad $\mT$ on an object $A$ is \bemph{unique to 
within a unique isomorphism}. In our case it means that, given two
action $\freccia{\mT A}{a,a'}{A}$ there is a unique invertible 2-cell
$\duefreccia{a}{\alpha}{a'}$ such that
$\freccia{(A,a)}{(1_A,\alpha)}{(A,a')}$ is a morphism of
$\mT$-algebras (in particular it is an isomorphism of
$\mT$-algebras). In this case we will say that the
$\mT$-\bemph{algebra structure is essentially unique}.  

More precisely a 2-monad $(\mT,\mu, \eta)$ is said
\bemph{property-like}, if it satisfies the following conditions: 
\begin{itemize}
\item for every $\mT$-algebra $(A,a)$ and $(B,b)$, and for every
invertible 1-cell $\freccia{A}{f}{B}$ there exists a unique
invertible 2-cell $\ovln{f}$ 
\[\xymatrix@+1pc{
\mT A \ar[d]_{a} \ar[r]^{\mT f} \xtwocell[r]{}<>{_<5>\ovln{f}}& \mT B \ar[d]^b\\
A \ar[r]_f &B
}\]
such that $\freccia{(A,a)}{(f,\ovln{f})}{(B,b)}$ is a morphism of
$\mT$-algebras;
\item  for every $\mT$-algebra $(A,a)$ and $(B,b)$, and for every
  1-cell $\freccia{A}{f}{B}$ if there exists a 2-cell $\ovln{f}$ 
\[\xymatrix@+1pc{
\mT A \ar[d]_{a} \ar[r]^{\mT f} \xtwocell[r]{}<>{_<5>\ovln{f}}& \mT B \ar[d]^b\\
A \ar[r]_f &B
}\]
such that $\freccia{(A,a)}{(f,\ovln{f})}{(B,b)}$ is a lax morphism of $\mT$-algebras, then it is the unique 2-cell with such property.
\end{itemize}
We conclude this section recalling a stronger property on a 2-monads $(\mT,\mu,\eta)$ on $\mA$ which implies that $\mT$ is property-like: a 2-monad $(\mT,\mu, \eta)$ is said \bemph{lax-idempotent}, if for every $\mT$-algebras $(A,a)$ and $(B,b)$, and for every 1-cell $\freccia{A}{f}{B}$ there exists a unique 2-cell $\ovln{f}$
\[\xymatrix@+1pc{
\mT A \ar[d]_{a} \ar[r]^{\mT f} \xtwocell[r]{}<>{_<5>\ovln{f}}& \mT B \ar[d]^b\\
A \ar[r]_f &B
}\]
such that $\freccia{(A,a)}{(f,\ovln{f})}{(B,b)}$ is a lax morphism of $T$-algebras.
In particular every lax-idempotent monad is property like. See \citep[Proposition 6.1]{OPLS}.

\section{Primary and existential doctrines}
The notion of hyperdoctrine was introduced by F.W. Lawvere in a series
of seminal papers \citep{AF,EHCSAF}. We recall from \citep{EQC}
some definitions which will be useful in the following. The reader can
find all the details about the theory of elementary and existential
doctrines also in \citep{EQC,QCFF,UEC}, and we refer to \citep{Frey2014AFS} for a detailed analysis of cocompletions of doctrines.

\begin{definition}
Let $\mC$ be a category with finite products. A \bemph{primary doctrine} is a functor
$\doctrine{\mC}{P}$ from the opposite of the category $\mC$ to the category of inf-semilattices.
\end{definition}

\begin{definition}\label{def elementary doctrine}
A primary doctrine $\doctrine{\mC}{P}$ is \bemph{elementary} if for every object $A$ in $\mC$ there exists an element $\delta_A$ in the fibre $P(A\times A)$ such that
\begin{enumerate}
\item the assignment
\[\Einv_{\angbr{\id_A}{\id_A}}(\alpha):=P_{\pr_1}(\alpha)\wedge \delta_A\]
for $\alpha$ in the fibre $P(A)$ determines a left adjoint to $\freccia{P(A\times A)}{P_{\angbr{\id_A}{\id_A}}}{P(A)}$;
\item for every morphism $e$ of the form $\freccia{X\times A}{\angbr{\pr_1,\pr_2}{\pr_2}}{X\times A\times A}$ in $\mC$, the assignment
\[ \Einv_{e}(\alpha):= P_{\angbr{\pr_1}{\pr_2}}(\alpha)\wedge P_{\angbr{\pr_2}{\pr_3}}(\delta_A)\]
for $\alpha$ in $P(X\times A)$ determines a left adjoint to $\freccia{P(X\times A \times A)}{P_e}{P(X\times A)}$.
\end{enumerate}
\end{definition}
\begin{definition}\label{def existential doctrine}
A primary doctrine $\doctrine{\mC}{P}$ is \bemph{existential} if, for every object $A_1$ and $A_2$ in $\mC$, for any projection $\freccia{A_1\times A_2}{{\pr_i}}{A_i}$, $i=1,2$, the functor
\[ \freccia{P(A_i)}{{P_{\pr_i}}}{P(A_1\times A_2)}\]
has a left adjoint $\Einv_{\pr_i}$, and these satisfy:
\begin{enumerate}
\item \bemph{Beck-Chevalley condition:} for any pullback diagram
\[
\quadratocomm{X'}{A'}{X}{A}{{\pr'}}{f'}{f}{{\pr}}
\]

with $\pr$ and $\pr'$ projections, for any $\beta$ in $P(X)$ the canonical arrow 
\[ \Einv_{\pr'}P_{f'}(\beta)\leq P_f \Einv_{\pr}(\beta)\]
is an isomorphism;
\item \bemph{Frobenius reciprocity:} for any projection $\freccia{X}{{\pr}}{A}$, for any object $\alpha$ in $P(A)$ and $\beta$ in $P(X)$, the canonical arrow
\[ \Einv_{\pr}(P_{\pr}(\alpha)\wedge \beta)\leq \alpha \wedge \Einv_{\pr}(\beta)\]
in $P(A)$ is an isomorphism.
\end{enumerate}
\end{definition}

\begin{remark}
In an existential elementary doctrine, for every map $\freccia{A}{f}{B}$ in $\mC$ the functor $P_f$ has a left adjoint $\Einv_f$ that can be computed as
\[ \Einv_{\pr_2}(P_{f\times {\id}_B}(\delta_B)\wedge P_{\pr_1}(\alpha))\]
for $\alpha$ in $P(A)$, where $\pr_1$ and $\pr_2$ are the projections from $A\times B$.
\end{remark}
Observe that primary doctrines, elementary doctrines, and existential doctrines have a 2-categorical structure given as follow. We refer to \citep{EQC,QCFF,UEC} for more details.

\begin{definition}\label{def 2-categroy PD}
The class of primary doctrines $\PD$ is a 2-category, where:
\begin{itemize}
\item \bemph{0-cells} are primary doctrines;
\item \bemph{1-cells} are pairs of the form $(F,b)$
\[\duemorfismo{\mC}{P}{F}{b}{\mD}{R}\]
such that $\freccia{\mC}{F}{\mD}$ is a functor preserving products, and $\freccia{P}{b}{R\circ F^{\op}}$ is a natural transformation such that the functor $\freccia{P(A)}{b_A}{RF(A)}$ preserves all the structure for every object $A$ in $\mC$, i.e. $b_A$ preserves finite meets;
\item \bemph{2-cells} $\duefreccia{(F,b)}{\theta}{(G,c)}$ are  natural transformations $\freccia{F}{\theta}{G}$ such that for every object $A$ in $\mC$ and for every $\alpha$ in $P(A)$, we have 
\[b_A(\alpha)\leq R_{\theta_A}(c_A(\alpha)).\] 
\end{itemize}
\end{definition}
Similarly we can define two 2-full 2-subcategories of $\PD$: the 2-category of existential doctrines $\ED$, and the 2-category of elementary doctrines $\ElD$. In these cases one should require that the 1-cells preserve the appropriate structures, in particular 1-cells of $\ED$ are those pairs $(F,b)$ such that $b$ preserves the left adjoints along projections. The 1-cells of $\ElD$ are those pairs $\freccia{P}{(F,b)}{R}$ such that for every object $A$ in $\mC$ we have
\[b_{A\times A}(\delta_A)=R_{\angbr{F \pr_1}{F\pr_2}}(\delta_{FA})\]
where $\delta_A=\Einv_{\Delta_A}(\top_A)$. See \citep{EQC,QCFF,UEC} for more details.
\begin{examples}\label{examples}
The following examples are discussed in \citep{AF}.
\begin{enumerate}
\item Let $\mC$ be a category with finite limits. The functor \[\doctrine{\mC}{{\Sub_{\mC}}}\]
assigns to an object $A$ in $\mC$ the poset $\Sub_{\mC}(A)$ of subobjects of $A$ in $\mC$ and, for an arrow $\frecciasopra{B}{f}{A}$ the morphism $\freccia{\Sub_{\mC}(A)}{\Sub_{\mC}(f)}{\Sub_{\mC}(B)}$ is given by pulling a subobject back along $f$. The fiber equalities are the diagonal arrows, so this is an elementary doctrine. Moreover it is existential if and only if the category $\mC$ is regular. See \citep{FSF}.
\item Consider a category $\mD$ with finite products and weak pullbacks: the doctrine is given by the functor of weak subobjects 
\[\doctrine{\mD}{{\Psi_{\mD}}}\]
where $\Psi_{\mD}(A)$ is the poset reflection of the slice category $\mD/A$, and for an arrow $\frecciasopra{B}{f}{A}$, the homomorphism $\freccia{\Psi_{\mD}(A)}{\Psi_{\mD}(f)}{\Psi_{\mD}(B)}$ is given by a weak pullback of an arrow $\frecciasopra{X}{g}{A}$ with $f$. This doctrine is existential, and the existential left adjoint are given by the post-composition.
\item Let $\mS$ be a theory in a first order language $\lang$. We define a primary doctrine 
\[\doctrine{\mC_{\mS}}{LT_{\mS}}\]
where $\mC_{\mS}$ is the category of lists of variables and term substitutions:
\begin{itemize}
\item \bemph{objects} of $\mC_{\mS}$ are finite lists of variables $\vec{x}:=(x_1,\dots,x_n)$, and we include the empty list $()$;
\item a \bemph{morphism} from $(x_1,\dots,x_n)$ to $(y_1,\dots,y_m)$ is a substitution $[t_1/y_1,\dots, t_m/y_m]$ where the terms $t_i$ are built in $\lang$ on the variable $x_1,\dots, x_n$;
\item the \bemph{composition} of two morphisms $\freccia{\vec{x}}{[\vec{t}/\vec{y}]}{\vec{y}}$ and $\freccia{\vec{y}}{[\vec{s}/\vec{z}]}{\vec{z}}$ is given by the substitution
\[ \freccia{\vec{x}}{[s_1[\vec{t}/\vec{y}]/z_k,\dots, s_k[\vec{t}/\vec{y}]/z_k]}{\vec{z}}.\]
\end{itemize}
The functor $\doctrine{\mC_{\mS}}{LT_{\mS}}$ sends $(x_1,\dots,x_n)$ in the class $LT_{\mS}(x_1,\dots,x_n)$  of all well formed formulas in the context $(x_1,\dots,x_n)$. We say that $\psi\leq \phi$ where $\phi,\psi\in LT_{\mS}(x_1,\dots,x_n)$ if $\psi\vdash_{\mS}\phi$, and then we quotient in the usual way to obtain a partial order on $LT_{\mS}(x_1,\dots,x_n)$. Given a morphism of $\mC_{\mS}$ 
\[\freccia{(x_1,\dots,x_n)}{[t_1/y_1,\dots,t_m/y_m]}{(y_1,\dots,y_m)}\]
the functor ${LT_{\mS}}_{[\vec{t}/\vec{y}]}$ acts as the substitution
${LT_{\mS}}_{[\vec{t}/\vec{y}]}(\psi(y_1,\dots,y_m))=\psi[\vec{t}/\vec{y}]$.

The doctrine $\doctrine{\mC_{\mS}}{LT_{\mS}}$ is elementary exactly when $\mS$ has an equality predicate. For all the detail we refer to \citep{QCFF}, and for the case of a many sorted first order theory we refer to \citep{CLP}.
\end{enumerate}

\end{examples}
\section{Existential completion}\label{section existential completion}
In this section we construct an existential doctrine $\doctrine{\mC}{\compex{P}}$, starting from a primary doctrine $\doctrine{\mC}{P}$. 

Let $\doctrine{\mC}{P}$ be a fixed primary doctrine for the rest of
the section, and let $\Lambda\subset \mC_1$ be a subset of morphisms closed
under pullbacks, compositions and such that it contains the identity
morphisms.

For every object $A$ of $\mC$ consider the following preorder:
\begin{itemize}
\item the objects are pairs $(\frecciasopra{B}{g\in \Lambda}{A},\; \alpha\in PB)$;
\item $(\frecciasopra{B}{h\in \Lambda}{A},\; \alpha\in PB)\leq (\frecciasopra{D}{f\in \Lambda}{A},\; \gamma\in PD)$ if there exists $\freccia{B}{w}{D}$ such that
\[\xymatrix@+1pc{
& B \ar[dl]_{w} \ar[d]^h\\
D \ar[r]_{f} &A
}\]
commutes and $\alpha\leq P_w(\gamma)$.
\end{itemize}
It is easy to see that the previous data give a preorder. Let
$\compex{P}(A)$ be the partial order obtained by identifying two
objects when
$$(\frecciasopra{B}{h\in \Lambda}{A},\; \alpha\in PB)\gtreqless (\frecciasopra{D}{f\in \Lambda}{A},\; \gamma\in PD)$$
in the usual way. With abuse of notation we denote the equivalence
class of an element in the same way.

Given a morphism $\freccia{A}{f}{B}$ in $\mC$, let $\compex{P}_f(\frecciasopra{C}{g\in \Lambda}{B},\; \beta\in PC)$ be the object 
\[(\frecciasopra{D}{\pbmorph{g}{f}}{A},\; P_{\pbmorph{f}{g}}(\beta)\in PD)\]
where  
\[\xymatrix@+1pc{
D\ar[d]_{\pbmorph{f}{g}} \pullbackcorner \ar[r]^{\pbmorph{g}{f}} &A\ar[d]^f\\
C \ar[r]_g & B
}\]
is a pullback because $g\in \Lambda$. 
\begin{proposition}
Let $\doctrine{\mC}{P}$ be a primary doctrine. Then $\doctrine{\mC}{\compex{P}}$ is a primary doctrine, in particular:
\begin{enumerate}
\item[(i)] for every object $A$ in $\mC$, $\compex{P}(A)$ is a inf-semilattice;
\item[(ii)] for every morphism $\freccia{A}{f}{B}$ in $\mC$, $\compex{P}_f$ is well-defined and it is an homomorphism of inf-semilattices.

\end{enumerate}
\end{proposition}
\proof
$(i)$ For  every $A$ we have the top element $(\frecciasopra{A}{\id_A}{A},\; \top_A)$. Consider two elements $(\frecciasopra{A_1}{h_1}{A},\; \alpha_1\in PA_1)$ and $(\frecciasopra{A_2}{h_2}{A},\; \alpha_2\in PA_2)$. In order to define the greatest lower bound of the two objects consider a pullback
\[\xymatrix@+1pc{
\pbp{A_1}{A_2}{A}\ar[d]_{\pbmorph{h_2}{h_1}} \pullbackcorner \ar[r]^{\pbmorph{h_1}{h_2}} &A_2\ar[d]^{h_2}\\
A_1 \ar[r]_{h_1} & A
}\]
which exists because $h_1\in \Lambda$ (and $h_2\in \Lambda$). We claim that 
\[(\frecciasopra{\pbp{A_1}{A_2}{A}}{\;\;\;\;{h_1\pbmorph{h_2}{h_1}}}{A}, P_{\pbmorph{h_2}{h_1}}(\alpha_1) \wedge P_{\pbmorph{h_1}{h_2}}(\alpha_2) )\] 
is such an infimum. 
It is easy to check that 
\[(\frecciasopra{\pbp{A_1}{A_2}{A}}{\;\;\;\;{h_1\pbmorph{h_2}{h_1}}}{A}, P_{\pbmorph{h_2}{h_1}}(\alpha_1) \wedge P_{\pbmorph{h_1}{h_2}}(\alpha_2) )\leq (\frecciasopra{A_i}{h_i}{A},\; \alpha_i\in PA_i)\]
 for $i=1,2$. Next consider  $(\frecciasopra{B}{g}{A},\; \beta\in PB)\leq (\frecciasopra{A_i}{h_i}{A},\; \alpha_i\in PA_i)$ for $i=1,2$ and $g=h_iw_i$. Then there is a morphism $\freccia{B}{w}{\pbp{A_1}{A_2}{A}}$ such that
\[\xymatrix@+1pc{
B\ar@{-->}[rd]^w \ar@/_1.3pc/[rdd]_{w_1} \ar@/^1.3pc/[rrd]^{w_2}  \\
&\pbp{A_1}{A_2}{A}\ar[d]_{\pbmorph{h_2}{h_1}} \pullbackcorner \ar[r]^{\pbmorph{h_1}{h_2}} &A_2\ar[d]^{h_2}\\
&A_1 \ar[r]_{h_1} & A
}\]
commutes and $\beta\leq P_{w_1}(\alpha_1)\wedge
P_{w_2}(\alpha_2)=P_w(P_{\pbmorph{h_2}{h_1}}(\alpha_1) \wedge
P_{\pbmorph{h_1}{h_2}}(\alpha_2))$. Observe that the infimum is well defined, since if, for example, we have
\[ (\frecciasopra{A_2}{h_2}{A},\; \alpha_2\in PA_2) \gtreqless (\frecciasopra{A_3}{h_3}{A},\; \alpha_3\in PA_3)\]
then there exist $\freccia{A_2}{w_3}{A_3}$ and $\freccia{A_3}{w_4}{A_2}$ such that $h_3w_3=h_2$, $\alpha_2\leq P_{w_3}(\alpha_3)$, $h_2w_4=h_3$ and $\alpha_3 \leq P_{w_4}(\alpha_2)$. Therefore there exists $\freccia{\pbp{A_1}{A_2}{A}}{w_5}{\pbp{A_1}{A_3}{A}}$ 
\[\xymatrix@+1pc{
\pbp{A_1}{A_2}{A}\ar@{-->}[rd]^{w_5} \ar@/_1.3pc/[rdd]_{\pbmorph{h_1}{h_2}} \ar@/^1.3pc/[rrd]^{w_3\pbmorph{h_2}{h_1}}  \\
&\pbp{A_1}{A_3}{A}\ar[d]_{\pbmorph{h_3}{h_1}} \pullbackcorner \ar[r]^{\pbmorph{h_1}{h_3}} &A_3\ar[d]^{h_3}\\
&A_1 \ar[r]_{h_1} & A
}\]
such that
\[P_{\pbmorph{h_2}{h_1}}(\alpha_1) \wedge P_{\pbmorph{h_1}{h_2}}(\alpha_2)\leq P_{w_5}(P_{\pbmorph{h_3}{h_1}}(\alpha_1) \wedge P_{\pbmorph{h_1}{h_3}}(\alpha_3)).\]
Then we can conclude that
\[ (\frecciasopra{\pbp{A_1}{A_2}{A}}{\;\;\;\;{h_1\pbmorph{h_2}{h_1}}}{A}, P_{\pbmorph{h_2}{h_1}}(\alpha_1) \wedge P_{\pbmorph{h_1}{h_2}}(\alpha_2) )\leq (\frecciasopra{\pbp{A_1}{A_3}{A}}{\;\;\;\;{h_1\pbmorph{h_3}{h_1}}}{A}, P_{\pbmorph{h_3}{h_1}}(\alpha_1) \wedge P_{\pbmorph{h_1}{h_3}}(\alpha_3) ).\]
Using the same argument one can prove that 
\[(\frecciasopra{\pbp{A_1}{A_3}{A}}{\;\;\;\;{h_1\pbmorph{h_3}{h_1}}}{A}, P_{\pbmorph{h_3}{h_1}}(\alpha_1) \wedge P_{\pbmorph{h_1}{h_3}}(\alpha_3) )\leq (\frecciasopra{\pbp{A_1}{A_2}{A}}{\;\;\;\;{h_1\pbmorph{h_2}{h_1}}}{A}, P_{\pbmorph{h_2}{h_1}}(\alpha_1) \wedge P_{\pbmorph{h_1}{h_2}}(\alpha_2) ).\]
Therefore we can conclude that the infimum is well-defined.

$(ii)$ We first prove that, for every morphism $\freccia{A}{f}{B}$, $\compex{P}_f$ is a morphism of preorders. By showing this, $\compex{P}_f$ will be a well-defined morphism of partial orders since we identify two elements $\ovln{\alpha}$ and $\ovln{\beta}$ of $\compex{P}(B)$ if $\ovln{\alpha} \gtreqless \ovln{\beta}$. Consider $(\frecciasopra{C_1}{g_1\in \Lambda}{B},\; \alpha_1\in PC_1)\leq (\frecciasopra{C_2}{g_2\in \Lambda}{B},\; \alpha_2\in PC_2)$ with $g_2w=g_1$ and $\alpha_1 \leq P_w(\alpha_2)$.
We want to prove that \[(\frecciasopra{D_1}{\pbmorph{g_1}{f}}{A},\; P_{\pbmorph{f}{g_1}}(\alpha_1)\in PD_1)\leq(\frecciasopra{D_2}{\pbmorph{g_2}{f}}{A},\; P_{\pbmorph{f}{g_2}}(\alpha_2)\in PD_1).\]
We can observe that $g_2w\pbmorph{f}{g_1}=g_1\pbmorph{f}{g_1}=f\pbmorph{g_1}{f}$. Then there exists a unique $\freccia{D_1}{\ovln{w}}{D_2}$ such that the following diagram commutes
\[\xymatrix@+1pc{
D_1 \ar@{-->}[rd]^{\ovln{w}} \ar@/_1.3pc/[rdd]_{w\pbmorph{f}{g_1}} \ar@/^1.3pc/[rrd]^{\pbmorph{g_1}{f}}  \\
&D_2 \ar[d]_{\pbmorph{f}{g_2}} \pullbackcorner \ar[r]^{\pbmorph{g_2}{f}} &A\ar[d]^{f}\\
&C_2 \ar[r]_{g_2} & B.
}\]
Moreover $P_{\ovln{w}}(P_{\pbmorph{f}{g_2}}(\alpha_2))=P_{\pbmorph{f}{g_1}}(P_w(\alpha_2))\geq P_{\pbmorph{f}{g_1}}(\alpha_1)$, and it is easy to see that $\compex{P}_f$ preserves top elements.
Finally it is straightforward to prove that $\compex{P}_f(\ovln{\alpha}\wedge\ovln{\beta})=\compex{P}_f(\ovln{\alpha})\wedge \compex{P}_f(\ovln{\beta})$. 
\endproof

\begin{proposition}\label{prop P_delta has left adjoint}
Given a morphism $\freccia{A}{f}{B}$ of $\Lambda$, let 
\[\compex{\Einv}_f(\frecciasopra{C}{h}{A},\alpha\in PC):=(\frecciasopra{C}{fh}{B},\alpha\in PC)\]
when $(\frecciasopra{C}{h}{A},\alpha\in PC)$ is in $\compex{P}(A)$.
Then $\compex{\Einv}_{f}$ is left adjoint to $\compex{P}_{f}$.
\end{proposition}
\proof
Let $\ovln{\alpha}:=(\frecciasopra{C_1}{g_1}{B},\; \alpha_1\in PC_1)$ and $\ovln{\beta}:=(\frecciasopra{D_2}{f_2}{A},\; \beta_2\in PD_2)$. 
Now we assume that $\ovln{\beta}\leq \compex{P}_{f}(\ovln{\alpha})$. This means that
\[\xymatrix@+1pc{
D_2 \ar[dr]^{f_2} \ar[d]_{w}\\
D_1 \pullbackcorner \ar[d]_{\pbmorph{f}{g_1}}\ar[r]^{\pbmorph{g_1}{f}}& A \ar[d]^{f}\\
C_1 \ar[r]_{g_1} & B
}\]
and $\beta_2\leq P_w(P_{\pbmorph{f}{g_1}}(\alpha_1))$. Then we have 
\[\xymatrix@+1pc{
 D_2 \ar[dr]^{f f_2} \ar[d]_{\pbmorph{f}{g_1}w}\\
C_1 \ar[r]_g &B
}\]
and $\beta_2 \leq P_{w\pbmorph{f}{g_1}}(\alpha_1)$. Then
$\compex{\Einv}_{f}(\ovln{\beta})\leq \ovln{\alpha}$.

Now assume $\compex{\Einv}_{f}(\ovln{\beta})\leq \ovln{\alpha}$
\[\xymatrix@+1pc{
 D_2 \ar[dr]^{f f_2} \ar[d]_{\ovln{w}}\\
C_1 \ar[r]_{g_1} &B
}\]
with $\beta_2 \leq P_{\ovln{w}}(\alpha_1)$.
Then there exists $\freccia{D_2}{w}{D_1}$ such that the following diagram commutes
\[\xymatrix@+1pc{
D_2 \ar@{-->}[rd]^{w} \ar@/_1.3pc/[rdd]_{\ovln{w}} \ar@/^1.3pc/[rrd]^{f_2}  \\
&D_1 \ar[d]_{\pbmorph{f}{g_1}} \pullbackcorner \ar[r]^{\pbmorph{g_1}{f}} &A\ar[d]^{f}\\
&C_1 \ar[r]_{g_1} & B
}\]
and $\beta_1 \leq P_{\ovln{w}}(\alpha_1)=P_w(P_{\pbmorph{f}{g_1}}(\alpha_1))$. Then we can conclude that $\ovln{\beta}\leq \compex{P}_{f}(\ovln{\alpha})$.
\endproof

\begin{theorem}\label{theorem generalized existential doctrine}
For every primary doctrine $\doctrine{\mC}{P}$, $\doctrine{\mC}{\compex{P}}$ satisfies:
\begin{itemize}
\item[(i)] \bemph{Beck-Chevalley Condition}: for any pullback
\[\xymatrix@+1pc{
X'\ar[d]_{f'} \pullbackcorner \ar[r]^{g'}& A'\ar[d]^{f}\\
X\ar[r]_{g} &A
}\]
with $g\in \Lambda$ (hence also $g'\in \Lambda$), for any $\ovln{\beta}\in \compex{P}(X)$ the following equality holds
\[\compex{\Einv}_{g'}\compex{P}_{f'}(\ovln{\beta})=\compex{P}_f\compex{\Einv}_{g}(\ovln{\beta}).\]
\item[(ii)] \bemph{Frobenius Reciprocity}: for every morphism $\freccia{X}{f}{A}$ of $\Lambda$, for every element $\ovln{\alpha}\in \compex{P}(A)$ and $\ovln{\beta}\in \compex{P}(X)$, the following equality holds
\[\compex{\Einv}_f(\compex{P}_f(\ovln{\alpha})\wedge \ovln{\beta})=\ovln{\alpha}\wedge \compex{\Einv}_f(\ovln{\beta}).\]
\end{itemize}
\end{theorem}
\proof
$(i)$ Consider the following pullback square
\[\xymatrix@+1pc{
X'\ar[d]_{f'} \pullbackcorner \ar[r]^{g'}& A'\ar[d]^{f}\\
X\ar[r]_{g} &A
}\]
where $g,g'\in \Lambda$, and let $\ovln{\beta}:=(\frecciasopra{C_1}{h_1}{X},\;\beta_1\in PC_1)\in \compex{P}(X)$. Consider the following diagram 
\[\xymatrix@+1pc{
D_1   \pullbackcorner \ar[r]^{\pbmorph{h_1}{f'}} \ar[d]_{\pbmorph{f'}{h_1}} &X'\ar[d]_{f'} \pullbackcorner \ar[r]^{g'}& A'\ar[d]^{f}\\
C_1 \ar[r]_{h_1}& X\ar[r]_{g} &A.
}\]
Since the two square are pullbacks, then the big square is a pullback, and then 
\[( \frecciasopra{D_1}{g'\pbmorph{h_1}{f'}}{A},P_{\pbmorph{f'}{h_1}}(\beta_1))=(\frecciasopra{D_1}{\pbmorph{gh_1}{f}}{A},P_{\pbmorph{f}{gh_1}}(\beta_1))\]
and these are by definition
\[\compex{\Einv}_{g'}\compex{P}_{f'}(\ovln{\beta})=\compex{P}_f\compex{\Einv}_{g}(\ovln{\beta}).\]
Therefore the Beck-Chevalley Condition is satisfied.

$(ii)$ Consider a morphism $\freccia{X}{f}{A}$ of $\Lambda$, an element $\ovln{\alpha}:=(\frecciasopra{C_1}{h_1}{A},\;\alpha_1\in PC_1)$ in $\compex{P}(A)$, and an element $\ovln{\beta}=(\frecciasopra{D_2}{h_2}{X},\;\beta_2\in PD_2)$ in $\compex{P}(X)$. Observe that the following diagram is a pullback 
\[\xymatrix@+1pc{
\pbp{D_1}{D_2}{X} \ar[d]_{\pbmorph{(\pbmorph{h_1}{f})}{h_2}}\pullbackcorner \ar[r]^{\pbmorph{h_2}{(\pbmorph{h_1}{f}})} &D_1\pullbackcorner \ar[d]_{\pbmorph{h_1}{f}} \ar[r]^{\pbmorph{f}{h_1}}& C_1 \ar[d]^{h_1}\\
D_2\ar[r]_{h_2}& X \ar[r]_f & A
}\]
and this means that
\[\compex{\Einv}_f(\compex{P}_f(\ovln{\alpha})\wedge \ovln{\beta})=\ovln{\alpha}\wedge \compex{\Einv}_f(\ovln{\beta}).\]
Therefore the Frobenius Reciprocity is satisfied.
\endproof

\begin{remark}
Observe that Proposition \ref{prop P_delta has left adjoint} and Theorem \ref{theorem generalized existential doctrine} just rely on the closure of the class $\Lambda$ of morphisms under composition and pullback and on the values of functors in meet semilattices,
while the finite product structure of $\mC$ is not used.
\end{remark}

We recall a useful lemma, which allows us to apply the previous construction on the class of projections, in order to obtain an existential doctrine in the sense of Definition \ref{def existential doctrine}.
\begin{lemma}\label{lemma proj are closed under pb}
Let $\mC$ be a category with finite products. Then the class of projections is closed under pullbacks, compositions and it contains identities.
\end{lemma}
\proof
It is direct to check that projections compose and that identities are projections. We show that this class is closed under pullbacks. Consider a projection $\freccia{A\times B}{\pr_A}{A}$ and an arbitrary morphism $\freccia{C}{f}{A}$ of $\mC$. It is direct to verify that the square
\[\comsquare{A\times B\times C}{C}{A\times B}{A}{\;\;\;\pr_C}{\angbr{f\pr_C}{\pr_B}}{f}{\pr_A}\]
commutes and it is a pullback.
\endproof

\begin{corollary}\label{corollary the doctrine compex is existential}
Let $\doctrine{\mC}{P}$ be a primary doctrine. If $\Lambda$ is the class of projections then the doctrine $\doctrine{\mC}{\compex{P}}$ is existential.
\end{corollary}

\begin{remark}\label{rem generalization of definitio of existential doctrine}
In the case that $\Lambda$ is the class of the projections, then from a primary doctrine $\doctrine{\mC}{P}$, we can construct an existential doctrine $\doctrine{\mC}{\compex{P}}$ in the sense of Definition \ref{def existential doctrine}.
Therefore the notion of existential doctrine can be generalized in the sense that an existential doctrine can be defined as a pair \[(\doctrine{\mC}{P},\Lambda)\]
where $\doctrine{\mC}{P}$ is a primary doctrine and $\Lambda$ is a class of morphisms of $\mC$ closed by pullbacks, composition and identities, which satisfies the conditions of Theorem \ref{theorem generalized existential doctrine}.
\end{remark}
\begin{remark}
Let $\freccia{\mC^{\op}}{P}{\pos}$ be a functor where $\pos$ is the category
of posets with top element. We can apply the existential completion
since we have not used the hypothesis that $PA$ has infimum in the
proofs; we have proved that if it has a infimum it is preserved by the
completion. In this case we must avoid to require Frobenius reciprocity.

Since a poset of the category $\pos$ has a top element, one has an
injection from the doctrine $\freccia{\mC}{P}{\pos}$ into
$\freccia{\mC}{\compex{P}}{\pos}$.
From a logical point of view, one can think of extending a theory
without existential quantification to one with that quantifier,
requiring that the theorems of the previous theory are preserved. 

We refer to \citep{Hofstra2010} for a general presentation of constructions which freely add quantification to a fibration, and their applications in logic.
\end{remark}
%\section{Construction of the 2-functor $\funE$}\label{section 2-fu}
In the rest of the section we assume that the morphisms of $\Lambda$ are all
the projections, since by Lemma \ref{lemma proj are closed under pb} this class is closed under pullbacks, compositions and it contains identities. 

We define a 2-functor $\freccia{\PD}{\funE}{\PED}$ from the 2-category of primary doctrines to the 2-category of existential doctrines, see Definition \ref{def 2-categroy PD}, which sends a primary doctrine $\doctrine{\mC}{P}$ to the existential
doctrine $\doctrine{\mC}{\compex{P}}$. For all the standard notions
about 2-category theory we refer to \citep{HCA1,HOHC}. 
\begin{proposition}\label{prop E is a due funct}
Consider the category $\PD(P,R)$. We define 
\[ \freccia{\PD(P,R)}{\funE_{P,R}}{\PED(\compex{P},\compex{R})}\]
as follow:
\begin{itemize}
\item for every 1-cell $(F,b)$, $\funE_{P,R}(F,b):=(F,\compex{b})$, where $\freccia{\compex{P}A}{\compex{b}_A}{\compex{R}FA}$ sends an object $(\frecciasopra{C}{g}{A},\; \alpha)$ in the object $(\frecciasopra{FC}{Fg}{FA},\;b_C(\alpha))$;
\item for every 2-cell $\duefreccia{(F,b)}{\theta}{(G,c)}$, $\funE_{P,R}\theta$ is essentially the same. 
\end{itemize}
With the previous assignment $\funE$ is a 2-functor.
\end{proposition}
\proof
We prove that $\freccia{\compex{P}}{(F,\compex{b})}{\compex{R}}$ is a
1-cell of $\PED(\compex{P},\compex{R})$. We first prove that for every
$A\in \mC$, $\compex{b}_A$ preserves the order.

If $(\frecciasopra{C_1}{g_1}{A},\; \alpha_1)\leq (\frecciasopra{C_2}{g_2}{A},\; \alpha_2)$, we have a morphism $\freccia{C_1}{w}{C_2}$ such that the following diagram commutes
\[\xymatrix@+1pc{
& C_1 \ar[dl]_w \ar[d]^{g_1}\\
C_2 \ar[r]_{g_2}& A
}\]
and $\alpha_1 \leq P_w(\alpha_2)$. Since $b$ is a natural
transformation, we have that $b_{C_1}P_w=R_{Fw}b_{C_2}$. Then we can
conclude that $(\frecciasopra{FC_1}{Fg_1}{FA},\;
b_{C_1}(\alpha_1))\leq (\frecciasopra{FC_2}{Fg_2}{FA},\;
b_{C_2}(\alpha_2))$ because $Fg_2Fw=Fg_1$ and
$b_{C_1}(\alpha_1) \leq b_{C_1}P_w(\alpha_2)=R_{Fw}(b_{C_2}\alpha_2)$. Moreover, since $F$ preserves products, we can
conclude that $\compex{b}_A$ preserves inf.

One can prove that
$\freccia{\compex{P}}{\compex{b}}{\compex{R}F^{\op}}$ is a natural
transformation using the facts that $F$ preserves products, which is needed to preserve projections. Moreover we can easily see that $\compex{b}$ preserves the left adjoints along
projections. Then $(F,\compex{b})$ is a 1-cell of $\PED$.

Now consider a 2-cell $\duefreccia{(F,b)}{\theta}{(G,c)}$, and let $\ovln{\alpha}=(\frecciasopra{C_1}{g_1}{A},\;\alpha_1)$ be an object of $\compex{P}(A)$. Then 
\[\compex{b}_A(\ovln{\alpha})=(\frecciasopra{FC_1}{Fg_1}{FA},\; b_{C_1}(\alpha_1))\]
and 
\[\compex{R}_{\theta_A}\compex{c}_A(\ovln{\alpha})=(\frecciasopra{D_1}{\pbmorph{Gg_1}{\theta_A}}{FA},\;R_{\pbmorph{\theta_A}{Gg_1}}c_{C_1}(\alpha_1))\]
where
\[\xymatrix@+1pc{
D_1 \pullbackcorner \ar[d]_{\pbmorph{\theta_A}{Gg_1}} \ar[r]^{\pbmorph{Gg_1}{\theta_A}} & FA\ar[d]^{\theta_A}\\
GC_1 \ar[r]_{Gg_1}& GA.
}\]
Now observe that since $\freccia{F}{\theta}{G}$ is a natural transformation, there exists a unique $\freccia{FC_1}{w}{D_1}$ such that the diagram

\[\xymatrix@+1pc{
FC_1\ar@{-->}[rd]^{w} \ar@/_1.3pc/[rdd]_{\theta_{C_1}} \ar@/^1.3pc/[rrd]^{Fg_1}\\
&D_1 \pullbackcorner \ar[d]_{\pbmorph{\theta_A}{Gg_1}} \ar[r]^{\pbmorph{Gg_1}{\theta_A}} & FA\ar[d]^{\theta_A}\\
&GC_1 \ar[r]_{Gg_1}& GA
}\]
commutes and then
$b_{C_1}(\alpha_1)\leq R_{\theta_{C_1}}c_{C_1}(\alpha_1)=R_wR_{\pbmorph{\theta_A}{Gg_1}}c_{C_1}(\alpha_1)$. Therefore we can conclude that
$\compex{b}_A(\ovln{\alpha})\leq \compex{R}_{\theta_A}\compex{c}_A(\ovln{\alpha})$,
and then $\freccia{F}{\theta}{G}$ can is a 2-cell
$\duefreccia{(F,\compex{b})}{\theta}{(G,\compex{c}),}$ and
$\funE_{P,R}(\theta  \gamma)=\funE_{P,R}(\theta) \funE_{P,R}(\gamma)$.

Finally one can prove that the following diagram commutes observing that for every $(F,b)\in \PD(P,R)$ and $(G,c)\in \PD(R,D)$, $(GF,\compex{c}\star\compex{b})=(GF,\compex{(c\star b)})$
\[\xymatrix@+1pc{
\PD(P,R)\times \PD(R,D)\ar[d]_{\funE_{PR}\times \funE_{RD}} \ar[rr]^{c_{PRD}} && \PD(P,D)\ar[d]^{\funE_{PD}}\\
\PED(\compex{P},\compex{R})\times \PED(\compex{R},\compex{D})\ar[rr]_{c_{\compex{P}\compex{R}\compex{D}}} && \PED(\compex{P},\compex{D})
}\]
where $c_{PRD}$ and $c_{\compex{P}\compex{R}\compex{D}}$ denote the composition functors of the 2-categories $\PD$ and $\ED$, and the same for the unit diagram.
Therefore we can conclude that $\funE$ is a 2-functor.
\endproof

Now we prove the 2-functor $\freccia{\PD}{\funE}{\PED}$ given by the assignment $\funE (P)=\compex{P}$ and by the functors $\funE_{P,R}$ defined in Proposition \ref{prop E is a due funct}, is left adjoint to the functor $\freccia{\PED}{\funU}{\PD}$ which \emph{forgets} the existential structure, i.e. it sends $P$ to $\funU(P)=P$.
\begin{proposition}\label{eta is 2-nat existential}
Let $\doctrine{\mC}{P}$ be a primary  doctrine. Then 
\[\freccia{P}{(\id_{\mC},\iota_P)}{\compex{P}}\]
where $\freccia{PA}{{\iota_P}_A}{\compex{P}A}$ sends $\alpha$ into $(\frecciasopra{A}{\id_A}{A},\;\alpha)$ is a 1-cell of primary doctrines. Moreover the assignment
\[\freccia{\id_{\PD}}{\eta}{\funU \funE}\]
where $\eta_P:=(\id_{\mC},\iota_P)$, is a 2-natural transformation.
\end{proposition}
\proof
It is easy to prove that $\freccia{PA}{{\iota_P}_A}{\compex{P}A}$ preserves all the structures.
For every morphism $\freccia{A}{f}{B}$ of $\mC$, it one can see that the following diagram commutes
\[\xymatrix@+1pc{
PB  \ar[d]_{{\iota_P}_B} \ar[r]^{P_f} & PA\ar[d]^{{\iota_P}_A}\\
\compex{P}B \ar[r]_{\compex{P}_f} & \compex{P}A.
}\]
Then we can conclude that $\freccia{P}{(\id_{\mC},\iota_P)}{\compex{P}}$ is a 1-cell of $\PD$ and it is a direct verification the proof the $\eta$ is a 2-natural transformation. 
\endproof

\begin{proposition}\label{proposition def epsilo counit}
Let $\doctrine{\mC}{P}$ be an existential doctrine. Then 
\[\freccia{\compex{P}}{(\id_{\mC},\zeta_P)}{P}\]
where $\freccia{\compex{P}A}{{\zeta_P}_A}{PA}$ sends $(\frecciasopra{C}{f}{A},\;\alpha)$ in $\Einv_f(\alpha)$ is a 1-cell of existential doctrines. Moreover the assignment 
\[\freccia{\funE\funU}{\varepsilon}{\id_{\ED}}\]
where $\varepsilon_P=(\id_{\mC},\zeta_P)$, is a 2-natural transformation.
\end{proposition}
\proof
Suppose $(\frecciasopra{C_1}{g_1}{A},\; \alpha_1)\leq (\frecciasopra{C_2}{g_2}{A},\; \alpha_2)$, with $\freccia{C_1}{w}{C_2}$, $g_2w=g_1$ and $\alpha_1\leq P_w(\alpha_2)$. Then by Beck-Chevalley
we have the equality
\[ \Einv_{\pbmorph{g_2}{g_1}} P_{\pbmorph{g_1}{g_2}}(\alpha_2)=P_{g_1}\Einv_{g_2}(\alpha_2)\]
and 
\[  \alpha_1\leq P_w(\alpha_2)\leq  P_wP_{g_2}\Einv_{g_2}(\alpha_2)= P_{g_1}\Einv_{g_2}(\alpha_2).\]
Then
\[ \Einv_{g_1}(\alpha_1)\leq \Einv_{g_2}(\alpha_2)\]
since $\Einv_{g_1}\dashv P_{g_1}$, and $\top_A=\zeta_A(\frecciasopra{A}{\id_A}{A},\;\top_A)$. Now we prove the naturality of $\zeta_P$. Let $\freccia{A}{f}{B}$ be a morphism of $\mC$. Then the following diagram commutes
\[\xymatrix@+1pc{
\compex{P}B \ar[d]_{\zeta_B} \ar[r]^{\compex{P}_f}& \compex{P}A \ar[d]^{\zeta_A}\\
PB \ar[r]_{P_f}& PA
}\]
because it corresponds to the Beck-Chevalley condition. Moreover it is easy to see that $\zeta_P$ preserves
left-adjoints. Then we an conclude that for every elementary
existential doctrine $\doctrine{\mC}{P}$, $\zeta_P$ is a 1-cell of
$\ED$.

The proof of the naturality of $\varepsilon$ is a routine verification. One must use the fact that we are working in $\ED$, and then for every 1-cell $(F,b)$, $b$ preserves left-adjoints along the projections.
\endproof

\begin{proposition}
For every primary doctrine $\doctrine{\mC}{P}$ we have 
\[\varepsilon_{\compex{P}} \circ \compex{\eta_P}=\id_{P}.\]
\end{proposition}

\proof
Consider the following diagram
\[\xymatrix@+1pc{
\mC^{op} \ar[rrd]^{\compex{P}}_{}="a" \ar[d]_{\id_{\mC}^{op}}\\
\mC^{op}\ar[d]_{\id_{\mC}^{op}} \ar[rr]_<<<<{\compex{(\compex{P})}}^{}="b"_{}="c"&& \infsl\\
\mC\ar[rru]_{\compex{P}}^{}="d"
\ar_{\compex{\iota}}  "a";"b"
\ar_{\zeta_{\compex{P}}}  "c";"d"}
\]
and let $(\frecciasopra{C}{g}{A},\; \alpha\in PA)$ be an element of $\compex{P}A$. Then 
\[ \compex{\iota_P}_A(\frecciasopra{C}{g}{A},\; \alpha\in PC)=(\frecciasopra{A}{\id_A}{A},\;(\frecciasopra{C}{g}{A},\; \alpha\in PC)\in \compex{P}A)\]
and
\[{\zeta_{\compex{P}}}_A(\frecciasopra{A}{\id_A}{A},\;(\frecciasopra{C}{g}{A},\; \alpha\in PC)\in \compex{P}A)=\compex{\Einv}_{\id_A}(\frecciasopra{C}{g}{A},\; \alpha\in PC).\]
By definition of $\compex{\Einv}$ we have
\[\compex{\Einv}_{\id_A}(\frecciasopra{C}{g}{A},\; \alpha\in PC)=(\frecciasopra{C}{g}{A},\; \alpha\in PC).\]
Then we can conclude that for every $\doctrine{\mC}{P}$, we have $\varepsilon_{\compex{P}}\circ \compex{\eta_P}=\id_{\compex{P}}$.
\endproof
\begin{corollary}\label{cor 1}
$\varepsilon \funE \circ \funE{\eta}=\id_{\funE}$.
\end{corollary}

\begin{theorem}\label{theorem E left adj U}
The 2-functor $\funE$ is 2-adjoint to the 2-functor $\funU$.
\end{theorem}
\proof
It is direct to verify that for every existential doctrine $\doctrine{\mC}{P}$ we have 
\[\varepsilon_P \circ \eta_P=\id_{P}\]
and then $\funU \varepsilon \circ \eta\funU=\id_{\funU}$. Therefore, by Corollary \ref{cor 1}, we can conclude that the 2-functor $\funE$ is 2-adjoint to the forgetful functor $\funU$, where $\eta$ is the unit of this 2-adjunction, and $\varepsilon$ is the counit.
\endproof

\section{The 2-monad $\mT_e$}\label{section 2-monad}
In this section we construct a 2-monad $\freccia{\PD}{\mT_e}{\PD}$,
and we prove that every existential doctrine can be seen as an algebra
for this 2-monad. Finally we show that the 2-monad $T_e$ is
lax-idempotent.

We define:
\begin{itemize}
\item $\freccia{\PD}{\mT_e}{\PD}$ the 2-functor $\mT=\forg \circ \funE$;
\item $\freccia{\id_{\PD}}{\eta}{\mT_e}$ is the 2-natural transformation defined in Proposition \ref{eta is 2-nat existential};
\item $\freccia{\mT_e^2}{\mu}{\mT_e}$ is the 2-natural transformation
$\mu=\forg\varepsilon \funE$.
\end{itemize}
%\end{definition}

\begin{proposition}
$\mT_e$ is a 2-monad.

\end{proposition}
\proof
One can easily check that the following diagrams commute
\[ \comsquare{\mT_e^3}{\mT_e^2}{\mT_e^2}{\mT_e}{\mu\mT_e}{T_e\mu}{\mu}{\mu}\]

\[\xymatrix{ 
\id_{\PD}\circ \mT_e\ar[r]^{\eta \mT_e} \ar[rd]_{\id}&\mT_e^2 \ar[d]^{\mu}& \mT_e\circ \id_{\PD}\ar[l]_{\mT_e\eta}\ar[ld]^{\id}\\
& \mT_e. 
}\] 
\endproof

\begin{proposition}\label{prop P existential implies T-algebgra}
Let $\doctrine{\mC}{P}$ be an existential doctrine. Then $(P,(\id_{\mC},\zeta_P))$ is a $\mT_e$-algebra, where $\freccia{\compex{P}}{\varepsilon_P=(\id_{\mC},\zeta_P)}{P}$ is the 1-cell of existential doctrines defined in Proposition \ref{proposition def epsilo counit}, i.e. $\freccia{\compex{P}A}{{\zeta_P}_A}{PA}$ sends $(\frecciasopra{C}{f}{A},\;\alpha)$ to $\Einv_f(\alpha)$. 
\end{proposition}

\proof
It is a direct verification.
\endproof

\begin{proposition}\label{prop stric algebra are existential}
Let $\doctrine{\mC}{P}$ be an primary doctrine, and let $(P,(F,a))$ be a $\mT_e$-algebra. Then $\doctrine{\mC}{P}$ is existential, $F=\id_{\mC}$ and $a=\zeta_P$.
\end{proposition}
\proof
By the identity axiom for $\mT_e$-algebras, we know that $F$ must be the identity functor, and $a_A\iota_A=\id_{PA}$.
\[\xymatrix@+1pc{
P\ar[rd]_{\id_P} \ar[r]^{\eta_P} & \compex{P} \ar[d]^{(F,a)}\\
& P.
}\]
For every morphism $\freccia{A}{f}{B}$ of $\mC$, where $f$ is a projection, we claim that 
\[\Einv_f(\alpha):=a_B \compex{\Einv}_f \iota_A(\alpha)\]
is left adjoint to $P_f$.
Let $\alpha\in PA$ and $\beta\in PB$, and suppose that $\alpha\leq P_f(\beta)$. Then we have that \[(\frecciasopra{A}{f}{B},\;\alpha)\leq (\frecciasopra{B}{\id_B}{B},\;\beta)\]
in $\compex{P}B$ and $(\frecciasopra{A}{f}{B},\;\alpha)=\compex{\Einv}_f(\frecciasopra{A}{\id_A}{A},\;\alpha)$. Therefore, by definition of $\iota$, we have
\[\compex{\Einv}_f\iota_A(\alpha)\leq \iota_B(\beta).\]
Hence
\[a_B\compex{\Einv}_f\iota_A(\alpha)\leq a_B\iota_B(\beta)=\beta.\]
Now suppose that $\Einv_f(\alpha)\leq \beta$. Then 
\[a_B(\frecciasopra{A}{f}{B},\;\alpha)\leq \beta\]
so
\[P_f a_B(\frecciasopra{A}{f}{B},\;\alpha)\leq P_f(\beta).\]
By the naturality of $a$, we have
\[P_f a_B(\frecciasopra{A}{f}{B},\;\alpha)=a_A\compex{P}_f(\frecciasopra{A}{f}{B},\;\alpha).\]
Now observe that $\iota_A(\alpha)=(\frecciasopra{A}{\id_A}{A},\; \alpha)\leq \compex{P}_f(\frecciasopra{A}{f}{B},\;\alpha)$. Therefore we have that
%\[\alpha= a_A\iota_A(\alpha)\leq P_f a_B(\frecciasopra{A}{f}{B},\;\alpha)\leq P_f(\beta).\]
\[\alpha\leq P_f(\beta)\]
follows from the unit law and the naturality of $a.$

Now we prove that Beck-Chevalley holds.
Consider the following pullback
\[\xymatrix@+1pc{
X' \ar[d]_{f'} \pullbackcorner\ar[r]^{g'} & A' \ar[d]^f\\
X \ar[r]_g & A
}\]
and $\alpha\in PX$. Then we have
\[\Einv_{g'}P_{f'}(\alpha)=a_{A'}\compex{\Einv}_{g'}\iota_{X'}(P_{f'}\alpha)=a_{A'}(\frecciasopra{X'}{g'}{A'},P_{f'}(\alpha)).\]
Observe that 
\[(\frecciasopra{X'}{g'}{A'},P_{f'}(\alpha))=\compex{P}_f(\frecciasopra{X}{g}{A},\;\alpha)\] and since $a$ is a natural transformation, we have
\[ a_{A'}\compex{P}_f(\frecciasopra{X}{g}{A},\;\alpha)=P_fa_A(\frecciasopra{X}{g}{A},\;\alpha).\]
Finally we can conclude that Beck-Chevalley holds because
\[P_f \Einv_g(\alpha)=P_f a_A\compex{\Einv}_g \iota_X(\alpha)=P_f a_A(\frecciasopra{X}{g}{A},\;\alpha).\]
Hence
\[\Einv_{g'}P_{f'}(\alpha)=P_f\Einv_g(\alpha).\]
Now consider a projection $\freccia{A}{f}{B}$, and two elements $\beta\in PB$ and $\alpha\in PA$.
We want to prove that the Frobenius reciprocity holds.
\[ \Einv_f(P_f(\beta)\wedge \alpha)=a_B\compex{\Einv}_f\iota_A(P_f(\beta)\wedge \alpha)=a_B(\frecciasopra{A}{f}{B},\;P_f(\beta)\wedge \alpha)\]
and 
\[ \beta \wedge \Einv_f(\alpha)= a_B \iota_B (\beta) \wedge a_B(\frecciasopra{A}{f}{B},\; \alpha)\]
and
\[a_B \iota_B (\beta) \wedge a_B(\frecciasopra{A}{f}{B},\; \alpha)=a_B((\frecciasopra{B}{\id_B}{B},\; \beta)\wedge (\frecciasopra{A}{f}{B},\; \alpha)).\]
We can observe that
\[a_B((\frecciasopra{B}{\id_B}{B},\; \beta)\wedge (\frecciasopra{A}{f}{B},\; \alpha))=a_B(\frecciasopra{A}{f}{B},\; P_f(\beta)\wedge \alpha)\]
and conclude that 
\[ \Einv_f(P_f(\beta)\wedge \alpha)= \beta\wedge \Einv_f(\alpha).\]
Therefore the primary doctrine $\doctrine{\mC}{P}$ is existential.
Finally we can observe that 
\[a_A(\frecciasopra{C}{g}{A},\alpha)=a_A\compex{\Einv}_g(\frecciasopra{C}{\id_C}{C},\;\alpha)=a_A\compex{\Einv}_g \iota_C(\alpha)=\Einv_g(\alpha).%%\qedhere
\]
Observe that all the previous calculations just depend on the naturality of $a$ and its unit law.
\endproof

\begin{proposition}\label{prop 1-cell of PED iff 1-cell of T-alg}
Let $(P,(\id_{\mC},\zeta_P))$ and $(R,(\id_{\mD},\zeta_R))$ be two $\mT_e$-algebras. Then every morphism $\freccia{(P,(\id_{\mC},\zeta_P))}{(F,b)}{(R,(\id_{\mD},\zeta_R))}$  of $\mT_e$-algebras is a 1-cell of $\PED$. Moreover every 1-cell of $\PED$ induces a morphism of $\mT_e$-algebras.
\end{proposition}
\proof
Let $\freccia{(P,(\id_{\mC},\zeta_P))}{(F,b)}{(R,(\id_{\mD},\zeta_R))}$  be a 1-cell of $\mT_e$-algebras. By definition of morphism of $\mT_e$-algebras, the following diagram commutes
\[\xymatrix@+1pc{
\compex{P}\ar[r]^{(F,\compex{b})} \ar[d]_{(\id_{\mC},\zeta_P)} & \compex{R}\ar[d]^{(\id_{\mD},\zeta_R)}\\
P \ar[r]_{(F,b)} & R.
}\]
Then for every object $(\frecciasopra{C}{g}{A},\;\alpha\in PC)$ of $\compex{P}A$ we have

\[\Einv^R_{Fg}b_C(\alpha)=b_A\Einv^P_g(\alpha)\] 
and this means that for every projection $\freccia{C}{g}{A}$ the following diagram commutes
\[\xymatrix@+1pc{
PC \ar[r]^{\Einv^P_g} \ar[d]_{b_C}& PA \ar[d]^{b_A}\\
RFC \ar[r]_{\Einv^R_{Fg}} & RFA.
}\]
Similarly one can prove that  every 1-cell of $\PED$ induces a morphism of $\mT_e$-algebras.
\endproof
\begin{corollary}\label{cor T-alg=PED}
We have the following isomorphism
\[ \alg{\mT_e}\cong \PED \] 
\end{corollary}
\proof
It follows from Proposition \ref{prop 1-cell of PED iff 1-cell of T-alg} and Proposition \ref{prop stric algebra are existential}. 
\endproof

Now we are going to prove that the 2-monad $\freccia{\PD}{\mT_e}{\PD}$ is lax-idempotent. This means that the 2-monad $\mT_e$ has both uniqueness of algebra structure and uniqueness of morphism structure, and then we can say that the \emph{existential} structure for a doctrine is a \emph{property} of the doctrine.

\begin{theorem}\label{prop T è lax-id}
Let $(P,(\id_{\mC},\zeta_P))$ and $(R,(\id_{\mD},\zeta_R))$ be $\mT_e$-algebras, and let $\freccia{P}{(F,b)}{R}$ be a 1-cell of $\PD$. Then $((F,b),\id_F)$ is a lax-morphism of algebras, and the 2-cell of primary doctrines $\duefreccia{(\id_{\mD},\zeta_R)(F,\compex{b})}{\id_F}{(F,b)(\id_{\mC},\zeta_P)}$ is the unique 2-cell making $((F,b),\id_F)$ a lax-morphism. Therefore the 2-monad $\freccia{\PD}{\mT_e}{\PD}$ is lax-idempotent.
\end{theorem}
\proof
Consider the following diagram where, following the notation of Proposition \ref{proposition def epsilo counit}, $\varepsilon_P=(\id_{\mC},\zeta_P)$ and $\varepsilon_R=(\id_{\mD},\zeta_R)$

\[\xymatrix@+1pc{
\compex{P} \ar[r]^{(F,\compex{b})} \ar[d]_{\varepsilon_P} \xtwocell[r]{}<>{_<5>\;\;\;\;\id_F}& \compex{R} \ar[d]^{\varepsilon_R}\\
P \ar[r]_{(F,b)} &R.
}\]
We must prove that for every object $A$ of $\mC$ and every
$(\frecciasopra{C}{f}{A},\;\alpha)$ in $\compex{P}A$ 
\[ \Einv_{Ff}^Rb_C(\alpha)\leq b_A\Einv_f^P(\alpha)\]
but the previous property holds if and only if
\[ b_C(\alpha)\leq R_{Ff}b_A\Einv_f^P(\alpha)=b_CP_f\Einv_f^P(\alpha)\]
and this holds since $\alpha\leq P_f\Einv_f^P(\alpha)$.

Finally it is easy to see that
$\duefreccia{\varepsilon_R(F,\compex{b})}{\id_F}{(F,b)\varepsilon_P}$
satisfies the coherence conditions for lax-$\mT_e$-morphisms.

Now suppose there exists another 2-cell $\duefreccia{\varepsilon_R(F,\compex{b})}{\theta}{(F,b)\varepsilon_P}$ such that $((F,b),\theta)$ is a lax-morphism
\[\xymatrix@+1pc{
\compex{P} \ar[r]^{(F,\compex{b})} \ar[d]_{\varepsilon_P} \xtwocell[r]{}<>{_<5>\;\theta
}& \compex{R} \ar[d]^{\varepsilon_R}\\
P \ar[r]_{(F,b)} &R.
}\]
Then it must satisfy the following condition
\[\xymatrix@+1pc{
P \ar[d]_{\eta_A}\ar[r]^{(F,b)} &R\ar[d]^{\eta_B}&&P\ar[dd]_{1_P}\ar[r]^{(F,b)}& R\ar[dd]^{1_B}\\
\compex{P} \ar[d]_{\varepsilon_P} \ar[r]^{(F,\compex{b})} \xtwocell[r]{}<>{_<5>\theta}&\compex{R} \ar[d]^{\varepsilon_R}\ar@{}[rr]|{=}&& \\
P \ar[r]_{(F,b)} & R&& P\ar[r]_{(F,b)} &R
}\]
and this means that $\theta=\id_F$.
\endproof

\begin{remark}
Observe that the family
$\duefreccia{\id_{\compex{\compex{P}}}}{\lambda_P}{\eta_{\compex{P}}
\mu_P}$ defined as $\lambda_P:=\id_{\mC}$ is a 2-cell in $\ED$. 

It is clearly a natural transformation. We must check that for every $\alpha\in \compex{(\compex{P})}A$
 \[\alpha\leq {\iota_{\compex{P}}}_A{\zeta_{\compex{P}}}_A(\alpha).\]
Let $\alpha:=(\frecciasopra{C}{g}{A},\;(\frecciasopra{D}{f}{C},\;\beta\in PD))$.
Then we have
\[{\iota_{\compex{P}}}_A{\zeta_{\compex{P}}}_A(\alpha)={\iota_{\compex{P}}}_A(\frecciasopra{D}{gf}{A},\;\beta \in PD)=(\frecciasopra{A}{\id_A}{A},(\frecciasopra{D}{gf}{A},\;\beta \in PD)).\]
Now we want to prove that 
\[(\frecciasopra{D}{f}{C},\;\beta\in PD) \leq \compex{P}_g(\frecciasopra{D}{gf}{A},\;\beta \in PD).\]
To see this inequality we can observe that the following diagram
commutes
\[\xymatrix@+1pc{
D_2 \ar@{-->}[rd]^{w} \ar@/_1.3pc/[rdd]_{\id_D} \ar@/^1.3pc/[rrrd]^{f}  \\
&L \ar[d]_{m_1}\pullbackcorner\ar[r]^{m_2} &H \ar[d]_{h_1} \pullbackcorner\ar[r]^{h_2} & C \ar[d]^g\\
&D \ar[r]_{f} &C\ar[r]_g & A
}\]
since every square is a pullback, hence $P_w(P_{m_1}(\beta))=\beta$. 

Moreover one can check that 2-cell $\freccia{\id_{T_e^2}}{\lambda}{\eta T_e \mu}$ is a modification. See \citep{HCA1} for the formal definition of modifications.

Finally, observe that the 2-cell $\mu$ is left adjoint to $\eta \mT_e$, where the unit of the adjunction is $\lambda$ and the counit is the identity. This result follows from the fact that for every $\doctrine{\mC}{P}$, the first component of the 1-cells $\mu_P$, $\eta \mT_e$ are the identity functor, and since $\lambda_P$ is the identity natural transformation, when we look at the conditions of adjoint 1-cell in the 2-category $\Cat$, we can observe that all the components are identities.
\end{remark}
\begin{remark}\label{prop ps-algebra are existential}
By Proposition \ref{prop stric algebra are existential} and Proposition \ref{prop stric algebra are existential} we have that a doctrine is existential if and only if it has a structure of $\mT_e$-algebra, but we can show that this also holds in the general setting of pseudo-algebras: if $\doctrine{\mC}{P}$ is a primary doctrine, and if $(P,(F,a))$ is a pseudo-$\mT_e$-algebra, then the doctrine $\doctrine{\mC}{P}$ is existential (the converse holds since strict algebras are a particular case of pseudo-algebras).

We refer to \citep{A2CC,PHDTT} for all the details about the formal definition of pseudo-algebras, and their properties.

If $(P,(F,a))$ is a pseudo-algebra, then there exists an invertible 2-cell
\[\xymatrix@+1pc{
P\xtwocell[rd]{}<>{_<-2>\;\; a_{\eta}}\ar[rd]_{\id_P}\ar[r]^{\eta_A} & \compex{P} \ar[d]^{(F,a)}\\
& P
}\]
and by definition, it is a natural transformation
$\freccia{F}{a_{\eta}}{\id_{\mC}}$, and for every $A\in\mC$ and
$\alpha\in PA$ we have $a_A\iota_A(\alpha)=P_{{a_{\eta}}_A}(\alpha)$.

Now consider a morphism $\freccia{A}{f}{B}$ in $\mC$ and $\alpha\in PA$. We define
\[\Einv_f(\alpha):=P_{{{a_\eta}_A}^{-1}}a_B\compex{\Einv}_f\iota_A(\alpha).\]
Using the same argument of Proposition \ref{prop stric algebra are existential} we can conclude that the primary doctrine $\doctrine{\mC}{P}$ is existential.
\end{remark}

\begin{example}\label{example existential completion 1}
Consider the fragment $\lang$ of first order intuitionistic logic with logical  symbols $\top$ and $\wedge$, and let $\lang_{\exists}$  be the fragment whose logical symbols are $\top$, $\wedge$ and $\exists$. Then we have that, following the notation used in Example \ref{examples}, the syntactic doctrine
\[\doctrine{\mC_{\lang_{\exists}}}{LT_{\lang_{\exists}}}\]
is isomorphic to the existential completion
\[\doctrine{\mC_{\lang}}{\compex{LT_{\lang}}}\]
of the primary doctrine $\doctrine{\mC_{\lang}}{LT_{\lang}}$.

Observe that we have this isomorphism because the operation of extending a language with the existential quantification is a free operation on the logic, so by the known equivalence between doctrines and logic given by the internal language, see for example \citep{CLP}, and since by Theorem \ref{theorem E left adj U} the existential completion is a free completion, these two doctrines must be isomorphic.

More specific categorical definitions of internal language are in \citep{MCBDTTCIPT,RCSILL}.
\end{example}

\section{Exact completion for elementary doctrines}
It is proved in \citep{UEC}  that there is a biadjunction  $\EED\rightarrow \excat$ between the 2-category of elementary existential doctrines and the 2-category of exact categories given by the composition of the
following 2-functors: the first is the left biadjoint to the inclusion
of $\CEED$ into $\EED$, see \citep[Theorem 3.1]{UEC}. The second is the
biequivalence between $\CEED$ and the 2-category $\LFS$ of categories
with finite limits and a proper stable factorization system, see
\citep{FSF}. The third is provided in \citep{NRRFS}, where it is proved
that the inclusion of the 2-category $\Reg$ of regular categories
(with exact functors) into $\LFS$ has a left biadjoint. The last
functor is the biadjoint to the forgetful functor from the 2-category
$\excat$ into $\Reg$, see \citep{REC}.

In this section we combine these results with the existential
completion for elementary doctrines, by proving the following result.
\begin{proposition}
The elementary structure is preserved by the existential completion, in the sense that if $\doctrine{\mC}{P}$ is
an elementary doctrine, then $\doctrine{\mC}{\compex{P}}$ is an
elementary existential doctrine.
\end{proposition}

Let $\doctrine{\mC}{P}$ be an elementary doctrine, and consider its existential completion $\doctrine{\mC}{\compex{P}}$. Given two objects $A$ and $C$ of $\mC$ we define
\[ \freccia{\compex{P}(A\times C)}{\compex{\Einv}_{\Delta_A\times \id_C}}{\compex{P}(A\times A\times C)}\]
on $\ovln{\alpha}=(\frecciasopra{A\times C\times D}{\;\;\;\pr}{A\times C},\; \alpha\in P(A\times C\times D))$ as
\[\compex{\Einv}_{\Delta_A\times \id_C}(\ovln{\alpha}):=(\frecciasopra{A\times A\times C\times D}{\;\;\;\pr }{A \times A\times C},\; \Einv_{\Delta_A\times \id_{C\times D}}(\alpha)\in P(A\times A \times C\times D)).\]

\begin{remark}
We can prove that $\compex{\Einv}_{\Delta_A\times \id_C}$ is a well defined functor for every $A$ and $C$: consider two elements of $\compex{P}(A\times C)$
\[\ovln{\alpha}=(\frecciasopra{A\times C\times D}{\;\;\;\pr}{A\times C},\; \alpha\in P(A\times C\times D))\]
and
\[\ovln{\beta}=(\frecciasopra{A\times C\times B}{\;\;\;\pr'}{A\times C},\; \beta\in P(A\times C\times B))\]
and suppose that $\ovln{\alpha}\leq \ovln{\beta}$. By definition there exists a morphism $\freccia{A\times C \times D}{f}{B}$ such that the following diagram commutes
\[\xymatrix@+1pc{
& A\times C\times D \ar[d]^{\pr_{A\times C}}\ar[dl]_{\angbr{\pr_{A\times C}}{f}}\\
A\times C \times B \ar[r]_{\pr'_{A\times C}} & A \times C
}\]
and $\alpha \leq P_{\angbr{\pr_{A\times C}}{f}}(\beta)$.
Since the doctrine $\doctrine{\mC}{P}$ is elementary we have
\[\beta\leq P_{\Delta_A \times \id_{C\times B}}\Einv_{\Delta_A \times \id_{C\times B}}(\beta)\]
and then
\[\alpha\leq P_{\angbr{\pr_{A\times C}}{f}}(P_{\Delta_A \times \id_{C\times B}}\Einv_{\Delta_A \times \id_{C\times B}}(\beta)).\]
Now observe that $(\Delta_A\times \id_{C\times B})({\angbr{\pr_{A\times C}}{f}})=({\angbr{\pr_{A\times A\times C}}{f\pr_{A\times C\times D}}})(\Delta_A \times \id_{C\times D})$, and this implies 
\[\alpha\leq P_{\Delta_A \times \id_{C\times D}}(P_{\angbr{\pr_{A\times A\times C}}{f\pr_{A\times C\times D}}}\Einv_{\Delta_A\times \id_{C\times B}}(\beta)).\]
Therefore we conclude 
\[\Einv_{\Delta_A\times \id_{C\times D}}(\alpha)\leq P_{\angbr{\pr_{A\times A\times C}}{f\pr_{A\times C\times D}}}\Einv_{\Delta_A\times \id_{C\times B}}(\beta).\]
It is easy to observe that the last inequality implies
\[\compex{\Einv}_{\Delta_A\times \id_{C}}(\ovln{\alpha})\leq \compex{\Einv}_{\Delta_A\times \id_{C}}(\ovln{\beta}).\] 

\end{remark}

\begin{proposition}
With the notation used before the functor 
\[\freccia{\compex{P}(A\times C)}{\compex{\Einv}_{\Delta_A\times \id_C}}{\compex{P}(A\times A \times C)}\]
is left adjoint to the functor
\[\freccia{\compex{P}(A\times A\times C)}{\compex{P}_{\Delta_A\times \id_C}}{\compex{P}(A\times C).}\]
\end{proposition}

\proof 
Consider an element $\ovln{\alpha}\in \compex{P}(A\times C)$,
\[\ovln{\alpha}:=(\frecciasopra{A\times C\times B}{\;\;\;\pr}{A\times C},\; \alpha\in P(A\times C\times B))\]
and an element $\ovln{\beta}\in \compex{P}(A\times A \times C)$,

\[\ovln{\beta}:=(\frecciasopra{A\times A\times C\times D}{\;\;\;\pr'}{A\times A \times C},\; \beta \in P(A\times A \times C\times D))\]
and suppose that
\[\compex{\Einv}_{\Delta_A\times \id_C}(\ovln{\alpha})\leq \ovln{\beta}\]
which means that there exists $\freccia{A\times A\times C\times B}{f}{D}$ 
\[\xymatrix@+1pc{
& A \times A \times C\times B\ar[d]^{\pr_{A\times A \times C}} \ar[dl]_{\angbr{\pr_{A\times A \times C}}{f}} \\
A\times A\times C\times D \ar[r]_{\pr_{A\times A \times C}} & A\times A \times C
}\]
such that $\Einv_{\Delta_A \times \id_{C\times B}}(\alpha)\leq P_{\angbr{\pr_{A\times A \times C}}{f}}(\beta)$.
Therefore we have
\[ \alpha \leq P_{\Delta_A\times \id_{C\times B}} P_{\angbr{\pr_{A\times A \times C}}{f}}(\beta) \]
and since 
\[(\angbr{\pr_{A\times A \times C}}{f})(\Delta_A \times \id_{C\times B})=(\Delta_A\times \id_{C\times D})\pr_{A\times C\times D}(\angbr{\pr_{A\times A \times C}}{f})(\Delta_A\times \id_{C\times B})\]
we can conclude that 
\[ \alpha\leq P_{\pr_{A\times C\times D}(\angbr{\pr_{A\times A \times C}}{f})(\Delta_A\times \id_{C\times B})}(P_{\Delta_A\times \id_{C\times D}}(\beta)).\]
Then
\[\ovln{\alpha}\leq \compex{P}_{\Delta_A\times \id_C}(\ovln{\beta})\]
because
\[\compex{P}_{\Delta_A\times \id_C}(\ovln{\beta})=(\frecciasopralunga{A\times C \times D}{\pr_{A\times C}}{A\times C},\; P_{\Delta_A \times \id_{C\times D}}(\beta)).\]
In the same way we can prove that $\ovln{\alpha}\leq \compex{P}_{\Delta_A\times \id_C}(\ovln{\beta})$ implies $\compex{\Einv}_{\Delta_A\times \id_C}(\ovln{\alpha})\leq \ovln{\beta}$.
\endproof

\begin{proposition}
Let $\compex{\delta}_A$ be $\compex{\Einv}_{\Delta_A}(\ovln{\top}_A)$. For every element $\ovln{\alpha}$ of the fibre $\compex{P}(A\times C)$ we have
\[ \compex{\Einv}_{ \Delta_A\times \id_C}(\ovln{\alpha})= \compex{P}_{\angbr{\pr_2}{\pr_3}}(\ovln{\alpha})\wedge \compex{P}_{\angbr{\pr_1}{\pr_2}}(\compex{\delta}_A)\]
where $\pr_i$, $i=1,2,3$, are the projections from $A\times A\times C$. In particular we have
\[ \compex{\Einv}_{ \Delta_A}(\ovln{\beta})= \compex{P}_{\pr_2}(\ovln{\beta})\wedge \compex{\delta}_A\]
for every element $\ovln{\beta}$ of the fibre $\compex{P}(A\times A)$.
\end{proposition}
 
\proof
Let $\ovln{\alpha}=(\frecciasopralunga{A\times C\times D}{\;\;\pr_{A\times C}}{A\times C},\; \alpha\in P(A\times C\times D))$ be an element of the fibre $\compex{P}(A\times C)$. Observe that $\compex{P}_{\angbr{\pr_2}{pr_3}}(\ovln{\alpha})$ is the element
$$\compex{P}_{\angbr{\pr_2}{pr_3}}(\ovln{\alpha})=(\frecciasopralunga{A\times A\times C\times D}{\;\;\pr_{A\times A\times C}}{A\times A\times C},\; P_{\angbr{pr_2',\pr_3'}{pr_4'}}(\alpha))$$
where $\freccia{A\times A\times C\times D}{\angbr{pr_2',\pr_3'}{pr_4'}}{A\times C\times D}$. Moreover we have that 
\[\compex{P}_{\angbr{\pr_1}{\pr_2}}(\compex{\delta}_A)=(\frecciasopra{A\times A\times C}{\id}{A\times A\times C},\; P_{\angbr{\pr_1}{\pr_2}}(\delta_A)).\]
Therefore $\compex{P}_{\angbr{\pr_2}{\pr_3}} (\ovln{\alpha})\wedge \compex{P}_{\angbr{\pr_1}{\pr_2}}(\compex{\delta}_A)$ is the element
\[(\frecciasopralunga{A\times A\times C\times D}{\;\;\pr_{A\times A\times C}}{A\times A\times C},\; P_{\angbr{pr_2',\pr_3'}{pr_4'}}(\alpha)\wedge P_{\angbr{\pr_1'}{\pr_2'}}(\delta_A)).\]
Note that $P_{\angbr{pr_2',\pr_3'}{pr_4'}}(\alpha)\wedge P_{\angbr{\pr_1'}{\pr_2'}}(\delta_A)=\Einv_{\Delta_A\times \id_{C\times D}}(\alpha)$ because the doctrine $P$ is elementary, so we can conclude that
\[ \compex{\Einv}_{ \Delta_A\times \id_C}(\ovln{\alpha})= \compex{P}_{\angbr{\pr_2}{\pr_3}}(\ovln{\alpha})\wedge \compex{P}_{\angbr{\pr_1}{\pr_2}}(\compex{\delta}_A).\]
\endproof

\begin{corollary}
For every elementary doctrine $\doctrine{\mC}{P}$, the existential completion $\doctrine{\mC}{\compex{P}}$ is elementary and existential.
\end{corollary}

\begin{example}\label{example 3}
Using the same argument of Example \ref{example existential completion 1}, one can prove that the syntactic doctrine
\[ \doctrine{\mC_{\lang_{=,\exists}}}{LT_{\lang_{=,\exists}}}\]
is the existential completion of the syntactic doctrine
\[ \doctrine{\mC_{\lang_{=}}}{LT_{\lang_{=}}}\]
where $\lang_{=,\exists}$ is the Regular fragment of first order intuitionistic logic, and $\lang_{=}$ is the Horn fragment.
\end{example}

We combine the existential completion for elementary doctrines with
the completions stated at the begin of this section, obtaining a
general version  of the exact completion described in
\citep{TECH,UEC}. We can summarise this operation with the following
diagram
\[ \xymatrix{
\ElD \ar[r]& \EED \ar[r] & \CEED \ar[r] & \LFS \ar[r] & \Reg \ar[r] &\excat.
}\]
It is proved in \textit{loc.cit.} that given an elementary
existential doctrine $\doctrine{\mC}{P}$ the completion
$\EED\rightarrow \excat$ produces an exact category denoted by
$\theory_P$ and this category is defined following the same idea used
to define a topos from a tripos. See \citep{TT,TTT}.

We conclude giving a complete description of the exact category
$\theory_{\compex{P}}$ obtained from an elementary doctrine
$\doctrine{\mC}{P}$.

Given an elementary doctrine $\doctrine{\mC}{P}$, consider the
category $\theory_{\compex{P}}$, called \bemph{exact completion of the
elementary doctrine} $P$, whose

\bemph{objects} are pair $(A,\rho)$ such that $\rho$ is in $P(A\times A\times C)$ for some $C$ and satisfies:
\begin{enumerate}
\item there exists a morphism $\freccia{A\times A \times C}{f}{C}$ such that 
\[\rho \leq P_{\angbr{\pr_2}{\pr_1,f}}(\rho)\] in $P(A\times A\times C)$ where 
$\freccia{A\times A\times C}{\pr_1,\pr_2}{A}$;
\item there exists a morphism $\freccia{A\times A\times A \times C}{g}{C}$ such that 
\[P_{\angbr{\pr_1,\pr_2}{\pr_4}}(\rho)\wedge P_{\angbr{\pr_2,\pr_3}{\pr_4}}(\rho)\leq P_{\angbr{\pr_1,\pr_3}{g}}(\rho)\] where
$\freccia{A\times A\times A\times C}{\pr_1,\pr_2,\pr_3}{A}$;
\end{enumerate}

\bemph{a morphism} $\freccia{(A,\rho)}{\phi}{(B,\sigma)}$, where $\rho\in P(A\times A\times C)$ and $\sigma\in P(B\times B\times D)$, is an object $\phi$ of $P(A\times B \times E)$ for some $E$ such that
\begin{enumerate}
\item there exists a morphism $\freccia{A\times B\times E}{\angbr{f_1}{f_2}}{C\times D}$ such that
\[ \phi \leq P_{\angbr{\pr_1,\pr_1}{f_1}}(\rho)\wedge P_{\angbr{\pr_2,\pr_2}{f_2}}(\sigma)\]
where the $\pr_i$'s are the projections from $A\times B\times E$;
\item there exists a morphism $\freccia{A\times A \times B\times C\times E}{h}{E}$ such that
\[P_{\angbr{\pr_1,\pr_2}{\pr_4}}(\rho)\wedge P_{\angbr{\pr_2}{\pr_3,\pr_5}}(\phi)\leq P_{\angbr{\pr_1}{\pr_3,h}}(\phi)\]
where the $\pr_i$'s are the projections from $A\times A\times B\times C\times E$;
\item there exists a morphism $\freccia{A\times B \times B\times D\times E}{k}{E}$ such that
\[P_{\angbr{\pr_2,\pr_3}{\pr_4}}(\sigma)\wedge P_{\angbr{\pr_1}{\pr_2,\pr_5}}(\phi)\leq P_{\angbr{\pr_1}{\pr_3,k}}(\phi)\]
where the $\pr_i$'s are the projections from $A\times B\times B\times D\times E$;
\item there exists a morphism $\freccia{A\times B \times B\times E}{l}{D}$ such that
\[ P_{\angbr{\pr_1,\pr_2}{\pr_4}}(\phi)\wedge P_{\angbr{\pr_1}{\pr_3,\pr_4}}(\phi)\leq P_{\angbr{\pr_2}{\pr_3,l}}(\sigma)\]
where the $\pr_i$'s are the projections from $A\times B\times B\times E$;
\item there exists a morphism $\freccia{A\times C}{\angbr{g_1}{g_2}}{B\times E}$ such that 
\[ P_{\angbr{\pr_1,\pr_1}{\pr_2}}(\rho)\leq P_{\angbr{\pr_1,g_1}{g_2}}(\phi)\]
where the $\pr_i$'s are the projections from $A\times C$.
\end{enumerate}
The composition of two morphisms is defined following the same
structure of the tripos to topos.

Observe that, in particular in point $5$ of the previous construction, the existential quantifiers disappear, because the usual last condition of the tripos-to-topos construction, see \citep{UEC,TTT}, which is the requirement $P_{\angbr{\pr_1}{\pr_1}}(\rho)\leq \Einv_{pr_2}(\phi)$, in the case $P$ is of the form $\compex{P}$, is equivalent to the condition $5$ of our previous construction because of the definition of the order in the fibre $\compex{P}(A)$.

Finally we conclude with the following theorem which generalized the exact completion for an elementary existential doctrine to an arbitrary elementary doctrine.
\begin{theorem}\label{theorem ED is biadjoint to Xct}
The 2-functor $\excat\rightarrow \ElD$ that takes an exact category to the elementary doctrine of its subobjects has a left biadjoint  which associates the exact category $\theory_{\compex{P}}$ to an elementary doctrine $\doctrine{\mC}{P}$.
\end{theorem}

\begin{example}
Combining Example \ref{example 3} and \citep[Theorem 4.7]{TECH}, we have that an instance of the previous construction is provided by the exact completion of existential m-variational doctrines $\mathbf{Ex}_{(LT_{\lang_{=,\exists}})_{cx}}$ defined in \citep{TECH}, which is isomorphic to the exact category $\theory_{\compex{(LT_{\lang_{=}})}}$.
\end{example}
Non-syntactic examples of existential completions and exact categories built from them are left to future work.

\end{document}